\documentclass[10pt,a4paper,oneside]{article}
\usepackage{enumerate,amsthm,amsmath,amsfonts,psfrag}
\usepackage {latexsym, verbatim}
\usepackage {graphics,graphicx}
\usepackage{color}
\usepackage {fancyhdr}
\usepackage[usenames,dvipsnames]{xcolor}
\usepackage [width=17cm]{geometry}
\usepackage{caption, subfigure}
\usepackage[colorlinks=true, a4paper=true, pdfstartview=FitV,linkcolor=blue, citecolor=blue, urlcolor=blue]{hyperref}

\newcommand {\pd} {\partial}

\newcommand {\e} {\varepsilon}
\newcommand{\dd}{\,\mathrm{d}}
\newcommand{\sign}{\,\mathrm{sign}}
\newcommand{\lu}{\underline\lambda}
\newcommand{\lo}{\overline\lambda}

\newtheorem{lemma}{Lemma}[section]  

\newtheorem{definition}[lemma]{Definition}
\newtheorem{proposition}[lemma]{Proposition}
\newtheorem{remark}{Remark}[section]
\newtheorem{example}{Example}[section]

\title{Counter-propagating waves in a system of transport-reaction equations} 
\author{Angelika Manhart} 
\date{} 
\begin{document}

\maketitle

\begin{center}
Courant Institute of Mathematical Sciences, \\
New York University, 251 Mercer Street, \\
New York, NY 10012, USA\\
angelika.manhart@cims.nyu.edu

\vspace{0.5 cm}
{\em This is a manuscript draft - the final version will be updated on arXiv early 2018}
\end{center}
\vspace{0.5 cm}
\begin{abstract}
Hyperbolic transport-reaction equations are abundant in the description of movement of motile organisms. Here, we focus on system of four coupled transport-reaction equations that arises from an age-structuring of a species of turning individuals. The highlight consists of the explicit construction and characterization of counter-propagating traveling waves, patterns which have been observed in bacterial colonies. Stability analysis reveals conditions for the wave formation as well as pulsating-in-time spatially constant solutions.
\end{abstract}


\vskip 0.4cm


\section{Introduction}

Large-scale patterns such as aggregation, waves or other structures created by groups of animals or microbes have always been a popular study object for mathematical modelers and analysts. Traditionally most such models for the movement of organisms were of parabolic nature \cite{burger2007,edelstein1998,holmes1993,mogilner1999,keller1970}.  The last two decades have seen a rise in the use of kinetic and hyperbolic models to model biological systems \cite{bellomo2007,carrillo2008,hadeler1999,hillen2002,othmer1988} (bacterial movement \cite{Lutscher2002}, chemotaxis of D. discoideum \cite{hillen2000}), see also the review in \cite{Eftimie2012}. Kinetic models are often derived as limits of individual based models (IBMs, also often called agent-based or particle models). IBMs are particularly useful if detailed experimental data, such as the speeds, tracks and interactions of the moving agents are available, as very little approximations are necessary to translate the biology into a mathematical formalism. On the other hand IBMs can typically be analyzed by direct simulation alone, limiting their predictive and mechanistic insight (see e.g. discussion in \cite{Manhart2016,holcombe2012})\newline\par

Hyperbolic models are typically harder to analyze than their diffusion-containing counter-parts. One of the most famous of such models is the Goldstein-Kac model \cite{goldstein1951,kac1974} for a correlated random walk in 1D: Here a particles moves with speed $s$ to the right or the left and its change of direction is modeled as a Poisson process with constant intensity $l$. Denoting the density of right- and left-moving particles by $u(x,t)$ and $v(x,t)$ this yields
\begin{align}
\label{eqn:intro1}
&\pd_t u + s \pd_x u =lv-lu,\\
&\pd_t v - s \pd_x v =lu-lv.\nonumber
\end{align}
Here the reaction term is linear and it is well known that the system is equivalent to a damped wave equation for the sum of the two densities (applying the \textit{Kac-trick} \cite{kac1974}).\newline\par

Several works have described the big variety of spatio-temporal patterns that nonlinear, local and non-local variants of \eqref{eqn:intro1} can produce \cite{eftimie2007,Eftimie2012}, these include stationary and traveling pulses, zigzag pulses and traveling trains. In this work we focus on what has been termed ripples, a spatio-temporal pattern in which two families of densities form counter-propagating traveling waves. Each density constitutes a genuine traveling wave, i.e. wave collisions do not affect the wave profiles. Such waves have already been observed as piecewise constant functions in \cite{Lutscher2002}, where a nonlinear version of the Goldstein-Kac model was analyzed:
\begin{align}
\label{eqn:intro2}
&\pd_t u + s \pd_x u =l(u)v-l(v)u,\\
&\pd_t v - s \pd_x v =l(v)u-l(u)v.\nonumber
\end{align}
Here $l(\rho)$ is a density dependent turning function. This system has also been mentioned in \cite{Degond2016}. Even for such simple models as \eqref{eqn:intro2} several open questions exist, such as potential blow-up behavior for quadratic $l$, or characterizing the limiting travelling waves for sigmoidal $l$. We will elucidate the latter question in Sec.~\ref{sec:waves}.\newline\par
However the main part of this work deals with a generalization of \eqref{eqn:intro2} in which an age-structuring is introduced \cite{DegMan2017, hadeler1999,Eftimie2012}. Here age refers to any internal variable that changes in time and can potentially affect speeds and turning rates. In the example about myxobacteria given in \cite{DegMan2017} it was assumed that the time since the last reversal, the ``age'' a, influences the turning behavior:
\begin{align}
\label{eqn:intro3}
&\pd_t u + s \pd_x u+g \pd_a u =-l(V,a)u, \quad u(x,0,t)=\frac{1}{g}\int_0^\infty l(U,a)v\dd a, \quad U(x,t)=\int_0^\infty u(x,a,t)\dd a\\
&\pd_t v - s \pd_x v+g\pd_a v =-l(U,a)v, \quad v(x,0,t)=\frac{1}{g}\int_0^\infty l(V,a)u\dd a, \quad V(x,t)=\int_0^\infty v(x,a,t)\dd a\nonumber
\end{align}
where $g$ is the rate of aging. Allowing for only two age-groups and assuming that only the group of the higher age can actually reverse allows to eliminate $a$ as an independent variable and leads to four transport-reaction equations that are coupled through their nonlinear reaction terms (see Sec.~\ref{sec:model} below and \cite{DegMan2017} for the derivation). We will see that also this model allows for counter-propagating traveling waves, in this case with more complex profiles than piecewise constant. However, the waves come with a twist: While their shape stays the same, their composition in terms of the two age groups is modified by the oncoming wave.
\newline\par
The rest of the manuscript is organized as follows: In Sec.~\ref{sec:model} we present the model and connect it to existing models through asymptotic scaling, in Sec.~\ref{sec:no_space} we analyze the space-independent dynamics in terms of existence and stability of steady states. In particular we find a Hopf-bifurcation. We proceed in Sec.~\ref{sec:steady} to analyze the effct the transport term has on the stability of the steady states. In Sec.~\ref{sec:waves} we show the existence and construction of counter-propagating traveling waves. Finally in Sec.~\ref{sec:numerics} we compare the calculated wave profiles to numerical simulations.

\section{Model Introduction and Scaling}
\label{sec:model}

The following model is a generalization of that derived and presented in \cite{DegMan2017}. We describe individuals moving in a 1D interval, $x\in[0,L]$. We denote by $u_{0/1}(x,t)$ and $v_{0/1}(x,t)$ the densities of individuals moving with speed $s>0$ to the right and left respectively. Each group can exist in one of two states: a reversible state (subscript $1$) and a non-reversible state (subscript $0$), which fulfill
\begin{align}
\label{eqn:main}
& \pd_t u_0 + s\, \pd_x u_0=-g(v_0+v_1) u_0+l( u_0+u_1) v_1,\\
& \pd_t  u_1 + s\, \pd_x u_1=g(v_0+v_1) u_0-l(v_0+v_1) u_1,\nonumber\\
& \pd_t  v_0 - s \,\pd_x v_0=-g(u_0+u_1) v_0+l(v_0+v_1) u_1,\nonumber\\
& \pd_t  v_1 - s\, \pd_x v_1=g(u_0+u_1) v_0-l(u_0+u_1) v_1.\nonumber
\end{align}
The left-hand-sides describe the (linear) transport of the densities, the right-hand-side the reactions. Individuals cycle through the different groups; starting with a non-reversible, right-moving individual we have: $u_0\rightarrow u_1$ (particle ``ages'' into a reversible state), $u_1\rightarrow v_0$ (particle reverses, thereby changes direction and becomes non-reversible), $v_0\rightarrow v_1$ (particle ``ages'' into a reversible state) and finally $v_1\rightarrow u_0$. The function $l(\rho)>0$ describes a turning rate that depends on the total density of opposing individuals, as has been suggested in \cite{DegMan2017}. The aging function $g(\rho)>0$ can be interpreted as the reciprocal value of a refraction period during which the individuals are unable to change their direction. It's dependence on the density of oncoming individuals presents a novelty compared to \cite{DegMan2017}. The interpretation is that the internal ``clock'' that determines the reversability is affected by signals coming from the opposing individuals. In general other dependencies of $l$ and $g$ on the densities are possible, however the procedure to construct the traveling waves in Sec.~\ref{sec:waves} might break down then. In previous works \cite{DegMan2017, scheel2017} $l(\rho)$ was assumed to be sigmoidal. Here, we don't prescribe particular shapes of $l(\rho)$ and $g(\rho)$, but keep the sigmoidal shape as reference example in mind.
\newline\par\noindent 
We equip system \eqref{eqn:main} with periodic boundary conditions. As initial average mass we define $m_0>0$ to be
$$
2m_0 =\frac{1}{L}\int_0^L \left( u_0+ u_1+ v_0+ v_1\right)\dd x,
$$

\paragraph{Scaling.}
We use as reference time, space and density $L/s$, $L$ and $m_0$ respectively and define the dimensionless aging and reversal functions as
$$
\gamma(w)=\frac{L}{s}g(m_0\,w),\quad \lambda(w)=\frac{L}{s}l(m_0\,w).
$$
Using the same notation for the dimensionless densities, space and time, we now have: $x\in[0,1]$ and
$$
\int_0^1 (u_0+u_1+v_0+v_1)\dd x=2 ,
$$
The transformed system written in terms of $u=u_0+u_1$ and $v=v_0+v_1$, $u_1$ and $v_1$ reads
\begin{align}
\label{eqn:main_scaled}
& \pd_t u + \pd_x u=\lambda( u) v_1-\lambda(v)u_1,\\
& \pd_t  v - \pd_x  v=\lambda( v) u_1-\lambda(u)v_1,\nonumber\\
& \pd_t  u_1 + \pd_x u_1=\gamma(v) (u-u_1)-\lambda( v) u_1,\nonumber\\
& \pd_t  v_1 -  \pd_x  v_1= \gamma(u) (v-v_1)-\lambda( u) v_1.\nonumber
\end{align}

\noindent Throughout this work, we will often assume the following properties of $\lambda$ and $\gamma$:
\begin{align}
\label{eqn:lambdaGamma}
&\gamma, \lambda \in \mathcal{C}^k([0,\infty)),\quad k\geq 2,\quad  0<\underline\gamma\leq \gamma(\rho),\quad 0<\lu \leq \lambda(\rho),\quad \forall \rho\geq 0.
\end{align}

\paragraph{Limit of fast aging.} To link the model in \eqref{eqn:main_scaled} to other existing models, it is instructive to look at the limit of fast aging, i.e. we set $\gamma(w)\rightarrow \frac{1}{\e}\gamma(w)$ for a small parameter $\e$. If we develop the densities with respect to $\e$, i.e.
\begin{align*}
&u(x,t)=u^0(x,t)+\e u^1(x,t)+\mathcal{O}(\e^2),\quad v(x,t)=v^0(x,t)+\e v^1(x,t)+\mathcal{O}(\e^2),\\
&u_1(x,t)=u^0_1(x,t)+\e u_1^1(x,t)+\mathcal{O}(\e^2),\quad v_1(x,t)=v_1^0(x,t)+\e v_1^1(x,t)+\mathcal{O}(\e^2),
\end{align*}
we obtain the following sets of equations: \textit{To first order}
\begin{align}
\label{eqn:noAge}
& \pd_t u^0 + \pd_x u^0=\lambda( u^0) v^0-\lambda(v^0)u^0,\\
& \pd_t  v^0 - \pd_x  v^0=\lambda( v^0) u^0-\lambda(u^0)v^0,\nonumber\\
&u_1^0=u^0,\quad v_1^0=v^0,\nonumber
\end{align}
\textit{to second order}
\begin{align*}
& \pd_t u^1 + \pd_x u^1=\left[\lambda'(u^0)v^0-\lambda(v^0)\right]u^1-\left[\lambda'(v^0)u^0-\lambda(u^0)\right]v^1+\lambda(u^0)\lambda(v^0)(v^0-u^0),\\
& \pd_t  v^1 - \pd_x  v^1=\left[\lambda'(v^0)u^0-\lambda(u^0)\right]v^1-\left[\lambda'(u^0)v^0-\lambda(v^0)\right]u^1+\lambda(u^0)\lambda(v^0)(u^0-v^0),\\
&u_1^1=u^1-\frac{\lambda(u^0)}{\gamma(v^0)}v^0,\quad v_1^1=v^1-\frac{\lambda(v^0)}{\gamma(u^0)}u^0.
\end{align*}
The model without aging \eqref{eqn:noAge} (also termed \textit{memory-free model}) has been described already e.g. in \cite{DegMan2017,Lutscher2002}. In this model all individuals are always in a reversible state. Already in \cite{Lutscher2002} it has been observed that, for certain choices of $\lambda$, this model produces piecewise-constant counter-propagating waves. We will revisit \eqref{eqn:noAge} in Section \ref{sec:waves}, where we'll show some new properties.


\section{The space-independent dynamics}
\label{sec:no_space}

Before dealing with the effect of the transport operator, we study the space-independent system, which constitutes a system of four coupled nonlinear ODEs for $(u(t),v(t), u_1(t), v_1(t))$. Since $u(t)+v(t)\equiv 2$, we can rewrite the system in terms of only three equations, for $d(t)=(u(t)-v(t))/2$, $u_1(t)$ and $v_1(t)$.
\begin{align}
\label{eqn:main_ODE}
& \dot d=\lambda(1+d) v_1-\lambda(1-d)u_1,\\
& \dot u_1 =\gamma(1-d) (1+d-u_1)-\lambda( 1-d) u_1,\nonumber\\
& \dot v_1= \gamma(1+d) (1-d-v_1)-\lambda( 1+d) v_1.\nonumber
\end{align}
The next two lemmata collect some straight forward properties.


\begin{lemma}\label{lem:ODEexist} Let $\gamma, \lambda \in \mathcal{C}^k([0,\infty))$ $k\geq 2$ and assume that there exist constants $\lu, \underline\gamma>0$, such that 
$$\underline \gamma\leq \gamma(\rho),\quad \lu \leq \lambda(\rho),\quad \forall \rho\geq 0.$$
An (isotropic) steady state of \eqref{eqn:main_ODE} is given by
\begin{equation}
\label{eqn:SS1}
d=0,\quad u_1=v_1=\frac{\gamma(1)}{\gamma(1)+\lambda(1)}.
\end{equation}
A sufficient condition for the existence of an additional pair of steady states is given by
\begin{equation}
\label{eqn:SS2}
\tau:=\frac{\gamma(1)}{\lambda(1)}\left[\lambda'(1)-\lambda(1)\right]+\frac{\lambda(1)}{\gamma(1)}\left[\gamma'(1)-\gamma(1)\right]>0
\end{equation}
\end{lemma}
\noindent
\textbf{Proof.}
By solving the steady state equations of $u_1$ and $v_1$ for $u_1$ and $v_1$ respectively and substituting the expressions into the steady state equation for $d$, it is easy to see that finding steady states of \eqref{eqn:main_ODE} is equivalent to finding $\bar d\in[-1,1,]$ such that $\bar d$ is a root of
\begin{align*}
&G(d)=(1-d)Q_+(d)-(1+d)Q_-(d),\quad Q_\pm(d)=\frac{\lambda(1\pm d)\gamma(1\pm d)}{\lambda(1\pm d)+\gamma(1\pm d)}
\end{align*}
The remaining values are then given by
\begin{align}
\label{eqn:SS2_u1v1}
u_1=\frac{\gamma(1-\bar d)}{\gamma(1-\bar d)+\lambda(1-\bar d)}(1+\bar d), \quad v_1=\frac{\gamma(1+\bar d)}{\gamma(1+\bar d)+\lambda(1+\bar d)}(1-\bar d).
\end{align}
The fixed point associated to $G(0)=0$ is the isotropic steady state in \eqref{eqn:SS1}. Since $G$ is continuous, $G(d)=G(-d)$ and $G(1)<0$, a sufficient condition for the existence of two additional roots of $G$  is $G'(0)>0$, which is equivalent to $\tau>0$.
\qed


\begin{lemma}\label{lem:ODEstab} Let the assumptions of Lem.~\ref{lem:ODEexist} hold.\\
\begin{enumerate}
\item[a.] The (isotropic) steady state of \eqref{eqn:main_ODE} given by \eqref{eqn:SS1} is stable iff the following two conditions hold
\begin{align*}
&1.\quad  \tau<0\\
&2.\quad  \lambda'(1)<2\lambda(1),\,\text{or}\,\, \gamma(1)\not\in\left[\lambda'(1)-\lambda(1)-\sqrt{\lambda'(1)(\lambda'(1)-2\lambda(1))},\lambda'(1)-\lambda(1)+\sqrt{\lambda'(1)(\lambda'(1)-2\lambda(1))} \right], 
\end{align*}
where $\tau$ is defined in \eqref{eqn:SS2}.
\item[b.] A necessary condition for the stability of any non-isotropic steady state $\bar d\neq 0$, is $G'(\bar d)<0$.
\end{enumerate}
\end{lemma}
\noindent
\textbf{Proof.}
Part a: Linearizing around \eqref{eqn:SS1} yields a matrix whose characteristic polynomial is given by
\begin{align}
&p(z)=\left(z+l+g\right)\left[z^2+(l+g-2b)z+2(lg-bg-lc)\right]\\
&l=\lambda(1), \quad g=\gamma(1),\quad b=\lambda'(1)\frac{\gamma(1)}{\gamma(1)+\lambda(1)}, \quad c=\gamma'(1)\frac{\lambda(1)}{\gamma(1)+\lambda(1)}.
\end{align}
The eigenvalue $z=-l-g$ is always negative. For the remaining quadratic polynomials the Routh-Hurwitz criterion asserts that stability is equivalent to the positivity of both coefficients, which in turn is equivalent to the first claim of the lemma.\newline
Part b: Let $0 \neq \bar d\in[-1,1]$ be such that $G(\bar d)=0$, which is equivalent to
$$
\bar d=\frac{Q_+-Q_-}{Q_++Q_-}.
$$
Here and in the following we always evaluate all functions at $d=\bar d$. Using this, we can write
$$
G'(\bar d)=2\,\frac{Q_+'Q_--Q_-'Q_+}{Q_++Q_-}-(Q_++Q_-),
$$
where $'$ denotes the derivative w.r.t. to $d$. Linearizing around the non-isotropic steady state given by $\bar d$ and \eqref{eqn:SS2_u1v1} yields a matrix whose determinant, using the expression for $G'(\bar d)$ derived above, can be written as
$$
\left(\lambda(1+\bar d)+\gamma(1+\bar d)\right)\left(\lambda(1-\bar d)+\gamma(1-\bar d)\right)G'(\bar d),
$$
Since the determinant of a $3\times 3$ matrix is minus the constant term of the corresponding characteristic polynomial, the Routh-Hurwitz theorem requires it to be negative for stability. This finishes the proof of part b. \qed\newline


\begin{figure}[t!]
\centering
\subfigure[Bifurcation diagram]{\includegraphics[width=0.48\textwidth]{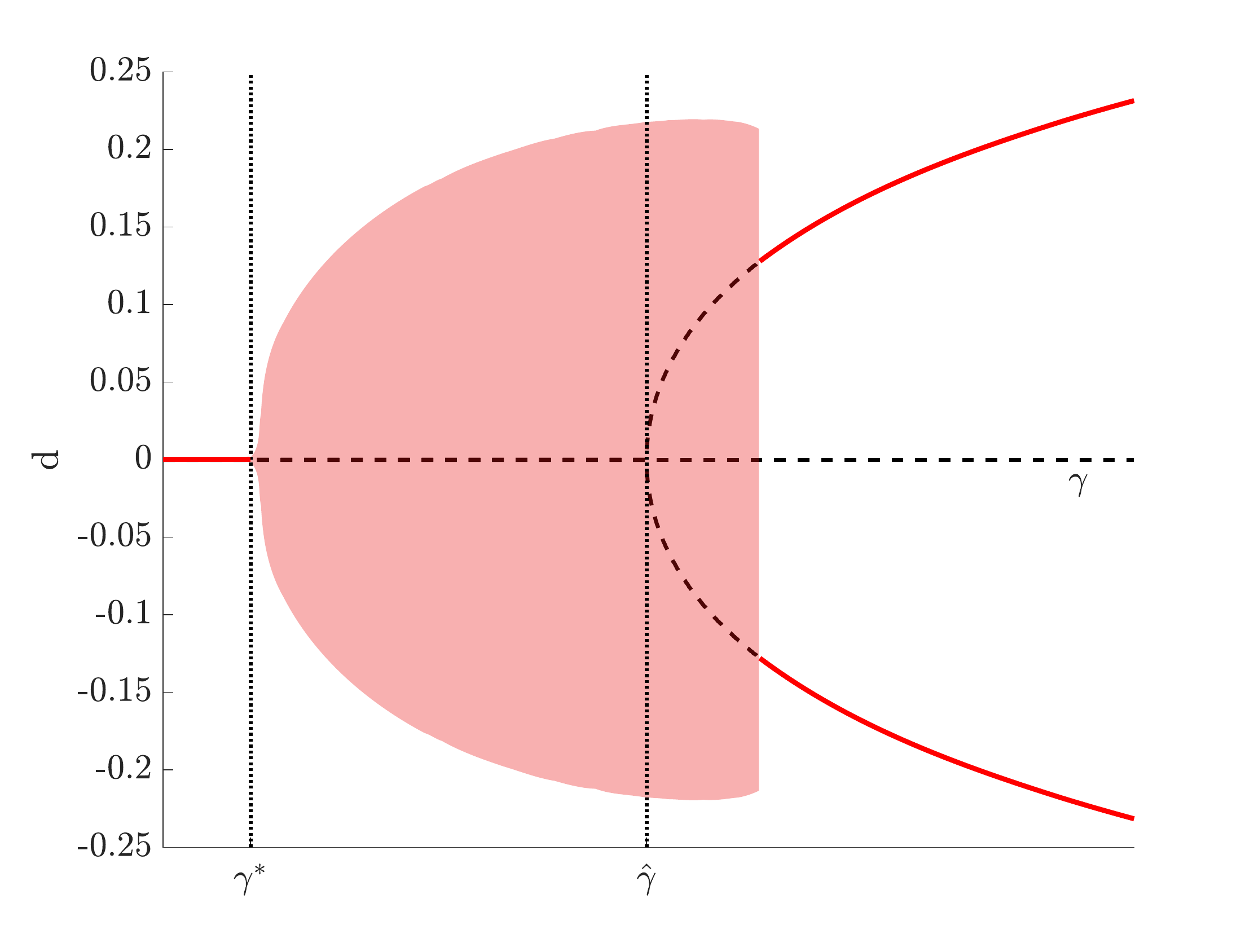}}
\subfigure[Limit sets in phase space]{\includegraphics[width=0.48\textwidth]{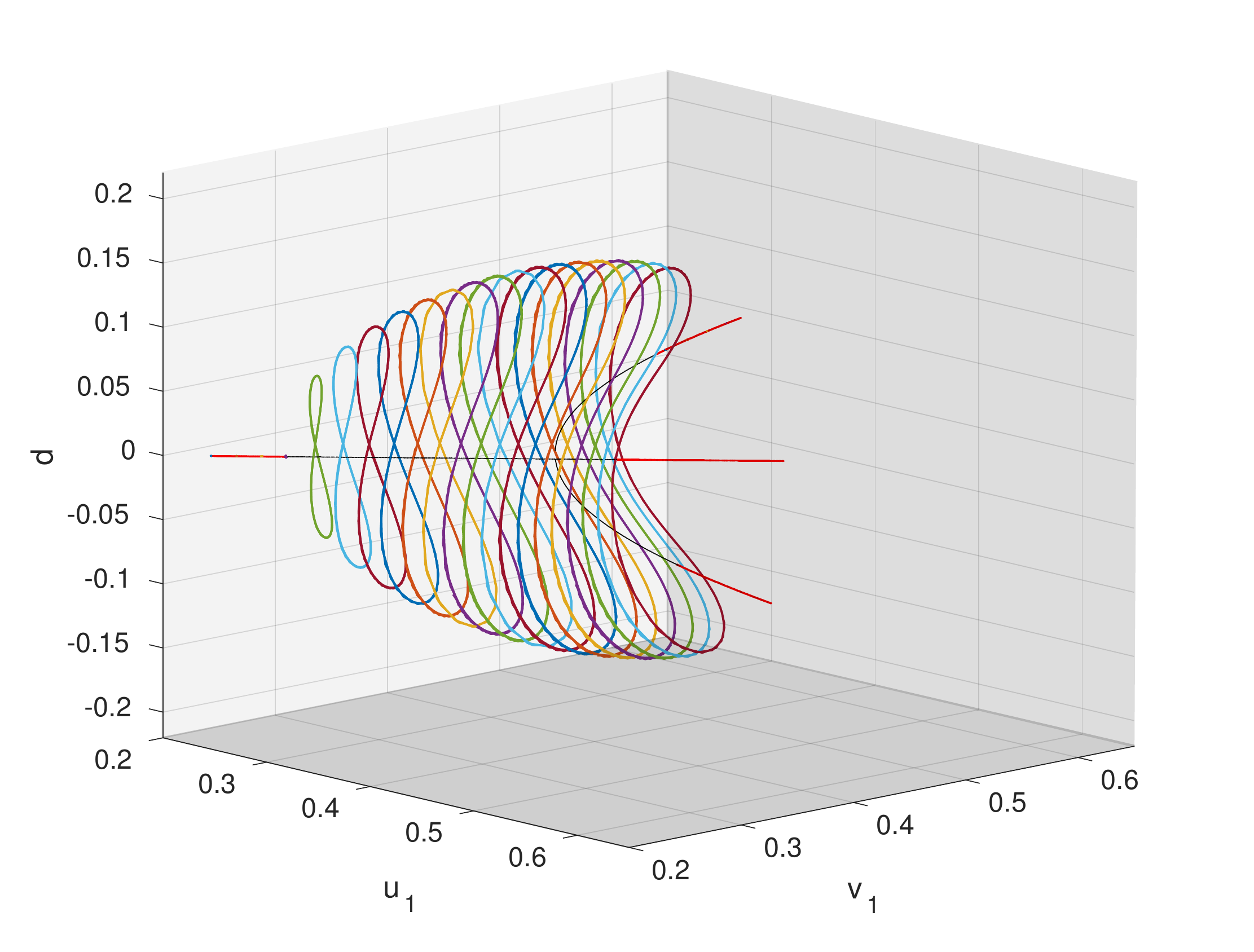}}
\caption{\textit{Hopf bifurcation.} (a) Shown is the bifurcation diagram for the parameter $\gamma$ using the limiting values of $d$ to represent the solution. Stable steady states are indicated by solid red lines, unstable steady states by dashed black lines. The region where a limit cycle is stable is shown in pale red. (b) The corresponding limit sets in $(u_1, v_1, d)$-space. Stable and unstable steady states are shown in thick red and thin black respectively, periodic limit cycles in many colors. Parameters are $(\lu, \lo,\alpha)=(2.5, 8,10)$, yielding $(\gamma^*, \hat\gamma, \gamma^{**})=(1.82, 3.24, 15.18)$. As initial conditions small perturbations from the isotropic steady state were used.} 
\label{fig:Hopf}
\end{figure}

It is noteworthy that the second condition in Lemma \ref{lem:ODEstab} allows for a situation where there exists only the isotropic steady state, but it is unstable. This indicates a more complicated long-term behavior, such as limit cycles. Let us given an

\begin{example}{\em Hopf bifurcation.}
\label{exp:Hopf} We use the following expressions for $\lambda(\rho)$ and $\gamma(\rho)$
\begin{align}
\label{eqn:lambda_hopf}
\lambda(\rho)=\lu+\frac{\lo-\lu}{1+\exp{(-\alpha(\rho-1))}}, \quad \gamma(\rho)\equiv\gamma
\end{align}
and interpret $\gamma>0$ as bifurcation parameter. $\lambda(\rho)$ is a sigmoid function taking values between $\lu$ and $\lo$. It has its steepest gradient, characterized by $\alpha$, at $\rho=1$. We define $\lambda_\pm=\frac{\lo\pm\lu}{2}>0$. To ensure that the first alternative of condition 2 in Lem.~\ref{lem:ODEstab}, part a is always violated we choose $\alpha>4\lambda_+/\lambda_-$. Expressing the remaining conditions in terms of $\gamma$ now yields that the isotropic steady state is stable iff the following two conditions hold
\begin{align*}
&1.\quad  \gamma<\frac{2\lambda_+^2}{\alpha\lambda_--2\lambda_+}=:\hat\gamma\\
&2.\quad \gamma\not\in \left[\frac{\alpha\lambda_--2\lambda_+-\sqrt{\alpha\lambda_-(\alpha\lambda_--4\lambda_+)}}{2} ,\frac{\alpha\lambda_--2\lambda_++\sqrt{\alpha\lambda_-(\alpha\lambda_--4\lambda_+)}}{2}\right]=:[\gamma^*, \gamma^{**}]
\end{align*}
Note that our choice of $\alpha$ also ensures that $\gamma^*<\hat\gamma<\gamma^{**}$. Finally we'd like to ensure that the sufficient condition for existence of a non-isotropic steady state of Lemma \ref{lem:ODEexist} is also necessary. An easy way to do this is to choose $\lambda_+>\lambda_-\sqrt{3}$ or equivalently $\lo<(2+\sqrt{3})\lu$. This can be checked easily, details are omitted. At the bifurcation point $\gamma=\gamma^*$, the isotropic fixed point becomes unstable and the pair of complex conjugated eigenvalues passes through the imaginary axis, i.e. we have a supercritical Hopf bifurcation. Fig. \ref{fig:Hopf} shows the bifurcation diagram and the limiting behavior of \eqref{eqn:main_ODE} obtained by varying $\gamma$. Biologically this indicates the existence of spatially constant, but pulsating in time solutions.
\end{example}


\section{Space-Homogeneous Steady States}
 \label{sec:steady}

\subsection{Stability of Space-Homogeneous Steady States}

A typical starting point of understanding pattern formation is to investigate the stability of space-homogeneous steady states. Their existence and values are discussed in Lem.~\ref{lem:ODEexist} and their stability in the space-indepedent system is discussed in Lem.~\ref{lem:ODEstab}. Next we discuss how the transport operator affects these stability criteria. 

\begin{proposition} 
\label{prop:stability} Let $\gamma, \lambda \in \mathcal{C}^k([0,\infty))$ $k\geq 2$ and assume that there exist constants $\lu, \underline\gamma>0$, such that 
$$\underline \gamma\leq \gamma(\rho),\quad \lu \leq \lambda(\rho),\quad \forall \rho\geq 0.$$
\begin{enumerate}
\item[a.] Let $(\bar u, \bar v, \bar u_1, \bar v_1)$ be the isotropic space-homogeneous steady state given by \eqref{eqn:SS1}. Assume that
\begin{equation}
\label{eqn:stability}
0<\lambda'(1)<\lambda(1),\quad \gamma'(1)<\gamma(1).
\end{equation}
Then $(\bar u, \bar v, \bar u_1, \bar v_1)$ is linearly stable. If $\lambda'(1)>\lambda(1)$, it is always unstable.
\item[b.] Let $\tau>0$ hold and let $(\bar u, \bar v, \bar u_1, \bar v_1)$ be the non-isotropic space-homogeneous steady state given by \eqref{eqn:SS2}. Necessary conditions for its linear stability are
\begin{align*}
\lambda'(\bar u)\bar v_1 < \frac{\gamma(\bar v) \lambda(\bar v)}{\gamma(\bar v)+\lambda(\bar v)}, \quad \lambda'(\bar v)\bar u_1 < \frac{\gamma(\bar u) \lambda(\bar u)}{\gamma(\bar u)+\lambda(\bar u)}, \quad G'(\bar u-\bar v)<0
\end{align*}
\end{enumerate}
\end{proposition}

\medskip 
\noindent 
{\bf Proof.}
We only show the proof of part a, part b uses the same techniques, but since it involves the use of the Routh Hurwitz criterion for complex polynomials, the calculations become somewhat cumbersome. We linearize around the first steady state with perturbations of zero total mass. Using matrix notation and calling the perturbations $r(x,t)\in \mathbb{R}^4$ we obtain the linear system
\begin{align*}
\pd_t r+T\pd_x r= Mr,\qquad 
&T=\begin{pmatrix}1 & 0& 0& 0\\ 0 & -1& 0& 0\\ 0 & 0& 1& 0\\ 0 & 0& 0& -1 \end{pmatrix},\quad M=\begin{pmatrix}b & -b & -l &l\\ -b & b & l & -l \\ g & c-b & -g-l & 0\\ c-b & g & 0& -g-l \end{pmatrix},
\end{align*}
where $l, g, b$ and $c$ are defined as in the proof of Lem.~\ref{lem:ODEstab}.
To asses the behavior of the system we look at perturbations of the form
\begin{align*}
r(x,t)=r_0\,e^{\xi t+i\,k\,x},\quad k=2\pi n,\quad n\in \mathbb{Z}, w_0\in\mathbb{R}^4\backslash \{0\},\qquad \text{yielding}\quad \left(M-T\,i\,k\right)r_0=\xi r_0.
\end{align*}
Expressing \eqref{eqn:stability} in terms of $l,g,b$ and $c$ gives
\begin{align}
\label{eqn:stability2}
0<b<\frac{g\,l}{g+l}, \quad c<\frac{g\,l}{g+l}
\end{align}
Hence we need to determine under what conditions $A=M-T\,i\,k$ has only eigenvalues $\xi$ with negative real part for all $k$.
We use the Routh-Hurwitz criterion for the (fourth order) characteristic polynomial of $A$, $p_A$. This yields the following conditions in terms of the the coefficients of $p_A(x)=x^4+a_3x^3+a_2 x^2+a_1 x+a_0$, where $a_i=a_i(g,l,b,c,k)$.
\begin{enumerate}
\item $a_i>0$
\item$ p=a_3a_2-a_1=p_1 k^2+p_0>0$
\item $q=a_3a_2a_1-a_1^2-a_3^2a_0=q_1k^2+q_0>0$
\end{enumerate}
where
\begin{alignat}{1}
\label{eqn:coeffs}
&a_0=k^2(k^2 + g^2+l^2+2l(b-c)), \quad a_1=2\left[k^2(g+l-b)+(g+l)(g\,l-b\,g-c\,l)\right],\\
&a_2=2k^2+(g+l)(g+l-2b)+2(g\,l-b\,g-c\,l), \quad a_3=2(g+l-b),\nonumber\\
&p_0=2(g+l-2b)\left[(g+l)(g+l-b)+g\,l-b\,g-c\,l\right],\quad  p_1=2(g+l-b),\nonumber\\
&q_0=4(g+l)(g\,l-b\,g-c\,l)(g+l-2b)\left[(g+l)(g+l-b)+g\,l-b\,g-c\,l)\right], \quad q_1=16(g+l-b)^2\left[g\,l-b(\,g+l)\right]\nonumber
\end{alignat}
It is easy to see that \eqref{eqn:stability2} implies
\begin{align}
\label{eqn:stability3}
g\,l-b\,g-c\,l>0,\quad g+l-2b>0,
\end{align}
and with the coefficients expressed as in \eqref{eqn:coeffs} this shows that \eqref{eqn:stability2} implies stability. On the other hand if $b>\frac{g\,l}{g+l}$, $q_1<0$ and hence we have instability for large $k$.
\qed

\begin{remark}
Comparing Prop.~\ref{prop:stability} to Lem.~\ref{lem:ODEstab} we see that the stability criteria for the ODE system are necessary, but not sufficient for the space-dependent system. A crucial difference is the additional requirement that
$$
\lambda'(1)<\lambda(1),
$$
which makes it possible to choose functions $\lambda(\rho)$ and $\gamma(\rho)$ and parameter regimes where the destabilization of the isotropic steady state requires the transport operator.
\end{remark}

\begin{example}{\em Constant aging function.} For density independent aging functions $\gamma(\rho)\equiv \gamma$, criteria $\lambda'(1)<\lambda(1)$ becomes both necessary and sufficient for the stability of the isotropic state. In such a case the reversal function $\lambda$ needs to have super-linear growth at 1 for the isotropic state to destabilize and the destabilization is independent of $\gamma$.  In \cite{DegMan2017}, where such a model was applied to myxobacterial movement, two wave forming criteria were determined numerically. The authors used a growing, piecewise quadratic, sigmoidal $\lambda(\rho)$, characterized by a minimal value $\lambda_m$, a maximal value $\lambda_M$ and an inflection density $\bar\rho$. Translating criterion \eqref{eqn:stability} back to dimensional variables  yields $\lambda'(m_0)m_0<\lambda(m_0)$ , and applying it to the model in \cite{DegMan2017} gives as destabilization conditions
$$
\lambda_M-3\lambda_m>0,\quad \frac{2\lambda_m}{\lambda_M-\lambda_m}< \frac{m_0}{\bar\rho}< \frac{2(\lambda_M-2\lambda_m)}{\lambda_M-\lambda_m}.
$$
Hence we recover the two wave formation criteria observed in \cite{DegMan2017}: The maximal reveral rate needs to be large enough compared to the minimal reversal rate and the average total density needs to be close to the inflection density.
\end{example}


\section{Counter-propagating Traveling Waves}
\label{sec:waves}

Upon simulation of the full system \eqref{eqn:main_scaled} (see also Sec.~\ref{sec:numerics}) a typical observed behavior is the formation of counter-propagating traveling waves for the two total densities $u(x,t)$ and $v(x,t)$. We therefore introduce the following definition. 

\begin{definition}
\label{def:travelingWaves} We call any solution $(u,v,u_1,v_1)$ of \eqref{eqn:main_scaled} a \textbf{Counter-propagating Traveling Wave Solution}  with speed $\pm c$, if there exist functions $P(\xi)$ and $M(\eta)$, such that
$$
u(x,t)=P(x-ct),\quad v(x,t)=M(x+ct).
$$
\end{definition}

This section is devoted showing the existence of such waves and explicitly constructing their shapes. For the memory-free system \eqref{eqn:noAge} we also present some new stability results.


\subsection{Waves in the Memory-free System}

 \begin{figure}[p]
\centering
\includegraphics[width=0.9\textwidth]{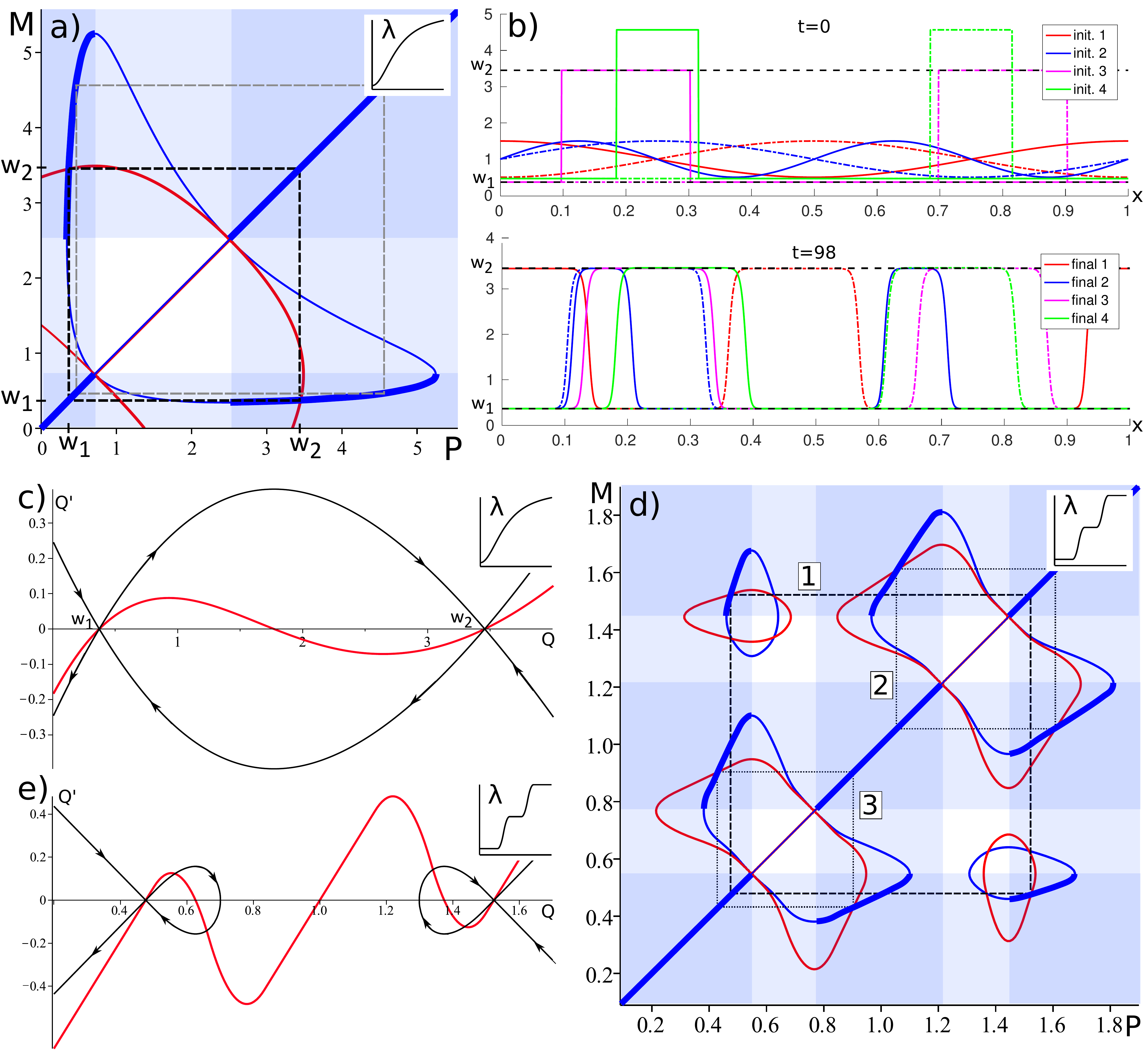}
\caption{\textit{Stable Traveling Waves in the Memory-free System.} a) Shown is the construction of stable candidate values of piecewise traveling waves in the memory-free system \eqref{eqn:noAge2} in $(u,v)=(P,M)$ space. The turning function used is described in Ex.~\ref{exp:sigmoid_ageless} and shown in the inset. Pairs of values that are roots of the right-hand-side in \eqref{eqn:noAge2} (condition {\bf A} in Rem.~\ref{rmk:ageless}) are shown in solid blue. Blue shading represents the linear stability condition  (condition {\bf B} in Rem.~\ref{rmk:ageless}), thick blue lines mark linearly stable pairs. Finally the integral condition {\bf C} is shown in red. The selected pair is depicted as thick dashed black square, a  linearly stable pair (not selected by the viscous limit) is shown in thin-dashed gray. b) Shown is the simulation outcome at time $t=98$ (lower row) for four different initial conditions (upper row). $u$ is shown as solid line, $v$ is shown as dash-dotted line. The constructed stable values $w_1$, $w_2$ from a) are shown as horizontal black dashed lines. Case 4 (green) corresponds to the pair shown in dashed gray in a). For this simulation we used $\Delta x=6.25\times 10^{-4}$ and a time step of $\Delta t=6.06\times 10^{-4}$ (compare Sec.~\ref{sec:numerics}). c) For the pair $(w_1,w_2)$ selected by the construction in a) the heteroclinic orbit in $(Q,Q')$-space is shown. This corresponds to the selected solution to the fast equation \eqref{eqn:fast}.   The inset shows the corresponding $\lambda(\rho)$. d) Construction of stable candidate values for Ex.~\ref{exp:doubleSigmoid}. The turning function $\lambda(\rho)$ used is shown in the inset. Linestyles and colors as in a), dashed black squares numbered 1, 2 and 3 represent three potentially stable pairs. e) As c), but for Ex.~\ref{exp:doubleSigmoid}. The turning function $\lambda(\rho)$ used is shown in the inset. The (homoclinic) orbit in $(Q, Q')$ space shown corresponds to the pair marked by 1 in d).}
\label{fig:ageless}
\end{figure}
 
 
Before analyzing the age-dependent system \eqref{eqn:main_scaled} we aim to understand the waves in the memory-free system \eqref{eqn:noAge}. By a slight abuse of notation we call the densities of right-moving and left-moving individuals again $u$ and $v$, i.e. we look at
\begin{align}
\label{eqn:noAge2}
& \pd_t u + \pd_x u=\lambda( u) v-\lambda(v)u,\\
& \pd_t  v - \pd_x  v=\lambda( v) u-\lambda(u)v.\nonumber
\end{align}
equipped with period boundary conditions on $[0,1]$ and $\int_0^1(u+v)\dd x =2$. This system has already been discussed in \cite{Lutscher2002, scheel2017}. Here we will summarize their findings and add a new result concerned with the stability of traveling waves.
\paragraph{Construction and linear stability.} We seek counter-propagating traveling waves in the sense of Def.~\ref{def:travelingWaves} of speed $\pm 1$. Making the ansatz $u(x,t)=P(x-t), v(x,t)=M(x+t)$ and defining as wave frames $\xi=x-t$ and $\eta=x+t$, one finds the condition
\begin{align}
\label{eqn:cond1}
\frac{P(\xi)}{\lambda(P(\xi))}=\frac{M(\eta)}{\lambda(M(\eta))}\equiv r>0.
\end{align}
Fig.~\ref{fig:ageless} shows all such solutions (blue curve) for a sigmoid $\lambda(\rho)$. For any given constant $r$, condition \eqref{eqn:cond1} will generically allow for only a finite set of values $P(\xi),M(\eta)\in\{w_1,\hdots w_K\}=:I_r$, hence any counter-propagating traveling wave will have to be piecewise constant. When assessing the linear stability of such waves, the simultaneous existence of two traveling waves frame provides a challenge and in general wouldn't allow for a standard perturbation ansatz. However, a work-around is to assume one density  to be completely constant (say $M(\eta)\equiv M$); now an ansatz in the frame of the other wave can be used, i.e.
$$
u(\xi,t)=P(\xi)+\e \tilde u(\xi,t),\quad v(\xi,t)=M(\xi)+\e \tilde v(\xi,t),
$$
for perturbations $\tilde u$ and $\tilde v$. Looking at the spectrum of the associated linear operator yields as stability condition (confer  \cite{scheel2017})
\begin{align}
\label{eqn:linearStabCond}
\lambda'(w_i)w_i<\lambda(w_i),
\end{align}
for any $w_i\in I_r$. An important observation is that $r$ parametrizes a whole family of stable tuples $I_r$ (compare Fig.~\ref{fig:ageless}a, thick blue lines), however upon simulation only a particular tuple associated to a particular value of $r$ seems to be selected. Note that also the (conserved) total mass cannot explain this selection since by shifting the jump points the mass can be varied without changing the wave heights.

\paragraph{Selection via diffusive regularization.} Already in \cite{DegMan2017} it has been noted that a small amount of density diffusion is necessary to be in agreement with the underlying particle model. We therefore add a small diffusion term to \eqref{eqn:noAge2} for $0<\e<<1$
\begin{align}
\label{eqn:noAge2_diffusion}
& \pd_t u + \pd_x u=\lambda( u) v-\lambda(v)u+\e^2\pd_x^2 u,\\
& \pd_t  v - \pd_x  v=\lambda( v) u-\lambda(u)v+\e^2\pd_x^2 v.\nonumber
\end{align}
Not we let $\e\rightarrow 0$ and consider the values selected by \eqref{eqn:noAge2_diffusion} to be the ``correct'' ones.  A similar procedure was also applied in \cite{anguige2009}. Let now $\{w_1,\hdots,w_K\}=I_r$ be a linearly stable solution tuple, i.e. fulfill \eqref{eqn:linearStabCond} and for a fixed $r>0$
\begin{align*}
\frac{w_i}{\lambda(w_i)}=r,\quad \forall i=1\hdots K
\end{align*}
We set $v(x,t)=M(x+t)\equiv w_1$ and search for solutions $u(x,t)=P(x-t)=P(\xi)$ of the first equation in \eqref{eqn:noAge2_diffusion}. We introduce a fast scale $\zeta=\xi/\e$ and define $Q(\zeta)=P(\xi)$. Substituting this into \eqref{eqn:noAge2_diffusion} yields
\begin{align}
\label{eqn:slow}
& 0=\lambda(P) w_1-\lambda(w_1)P+\e^2 P''\xrightarrow[]{\e\rightarrow 0} 0=\lambda(P) w_1-\lambda(w_1)P,\\
& 0=\lambda(Q) w_1-\lambda(w_1)Q+Q''.\label{eqn:fast}
\end{align}
We already know that for $\e\rightarrow 0$ on the slow scale $P$ will take values in $I_r$. Since we want to understand jumps between values in $I_r$, we search for solutions $Q$ that form a heteroclinic orbit connecting $w_1$ with any other $w_i\in I_r$, $i\neq 1$. Multiplying \eqref{eqn:fast} with $Q'$ and integrating with respect to $\zeta$ now yields
\begin{align}
C= \int^Q \left(\lambda(u)w_1-\lambda(w_1)u\right) \dd u+\frac{1}{2}(Q')^2,
\end{align}
for some constant $C$. Since we require $Q(\zeta)\rightarrow w_1$ and $Q(\zeta)\rightarrow w_i$ for $\zeta\rightarrow \pm \infty$, we obtain the integral condition
\begin{align}
\label{eqn:integral1}
 \int^{w_1} \left(\lambda(u)w_1-\lambda(w_1)u\right) \dd u= \int^{w_i} \left(\lambda(u)w_1-\lambda(w_1)u\right) \dd u.
\end{align}
This constitutes a necessary condition for ``correct'' tuples.

\begin{remark}{\em Alternative formulation and summary.}
\label{rmk:ageless} We can rewrite \eqref{eqn:integral1} as follows: By dividing the right-hand-side of \eqref{eqn:integral1} by $w_1$ and evaluating the second part of the integral we get
\begin{align*}
 \int^{w_i}\lambda(u)\dd u-\frac{\lambda(w_1)}{w_1}\frac{w_i^2}{2}.
\end{align*}
Since the fraction $\lambda(w_1)/w_1=\lambda(w_i)/w_i$ for all $w_i\in I_r$, we can replace it in the second term and obtain a formulation that is independent of $w_1$. We can now summarize the procedure to construct candidate values for stable (i.e. linearly stable and selected by the limit $\e\rightarrow 0$ in \eqref{eqn:noAge2_diffusion}) counter-propagating traveling waves as follows: Defining 
\begin{align}
\Lambda(\rho)=\frac{\rho}{\lambda(\rho)},\quad \Omega(\rho)=\int_0^\rho \lambda(u)\dd u-\frac{\lambda(\rho)\rho}{2},
\end{align}
any stable tuple $\{w_1,\hdots,w_K)$ has to fulfill 
\begin{align*}
&\textbf{A}:\,\, \Lambda(w_i)=\Lambda(w_j),\quad \textbf{B}:\,\, \Lambda'(w_i)>0, \quad  \textbf{C}:\,\, \Omega(w_i)=\Omega(w_j),\quad \forall i,j\in\{1,\hdots K\}.
\end{align*}
Note that {\bf{B}} is just a reformulation of \eqref{eqn:linearStabCond}. The sought solutions to \eqref{eqn:noAge2} are now of the form $u(x,t)=P(x-t)$, $v(x,t)=M(x+t)$ with $P(\xi),M(\eta)\in\{w_1,..w_K\}$. Fig.~\ref{fig:ageless} and Ex.~\ref{exp:sigmoid_ageless}, \ref{exp:doubleSigmoid} (below) illustrate this construction.
\end{remark}

The following lemma provides an easy way to find stable values for a class of turning functions $\lambda(\rho)$. In particular it applies to many sigmoid turning functions.


\begin{lemma}{\em Anti-symmetric turning function.} Let $\lambda(\rho)$ be anti-symmetric with respect to some density $\bar\rho>0$ in the following sense: There exist $0<\lu<\lo$ such that
\begin{align}
\label{eqn:symmetric}
\lo-\lambda(\bar\rho+\rho)=\lambda(\bar\rho-\rho)-\lu,\quad \forall \rho\in[\bar\rho,\bar\rho].
\end{align}
Then a stable pair of values $(w_1,w_2)$ can be found by replacing condition {\bf C} in Rem.~\ref{rmk:ageless} with
$$
w_2=2\bar\rho-w_1.
$$
\end{lemma}
\noindent
{\bf Proof.} Let $w_2=2\bar\rho-w_1$. We start by observing that condition {\bf A} together with \eqref{eqn:symmetric} implies that 
\begin{align}
\label{eqn:proof1}
\Lambda(w_1)=\Lambda(w_2)=\frac{2\bar\rho}{\lo+\lu}.
\end{align}
Performing a change of variables $u\rightarrow 2\bar\rho-u$ and using \eqref{eqn:symmetric} it is easy to see that
\begin{align*}
\hat\Omega(w)&:= \int_{\bar\rho}^{w} \left(\lambda(u)w_1-\lambda(w_1)u\right) \dd u=\hat\Omega(2\bar\rho-w)+(\bar\rho-w)\left[2\bar\rho\lambda(w_1)-(\lo+\lu)w_1\right]\\
&=\hat\Omega(2\bar\rho-w),
\end{align*}
where we have used \eqref{eqn:proof1} in the last line. Noting that condition {\bf C} is equivalent to $\hat\Omega(w_1)=\hat\Omega(w_2)$ finishes the proof.
\qed



\begin{example}{\em Linear and quadratic turning functions.}\label{exp:linear} If $\lambda(\rho)=a+b\rho$, i.e. linear, the effect is equivalent to a constant turning function $\lambda(\rho)=a$, since the terms involving $b$ cancel. In this case \eqref{eqn:noAge2} reduces to the Goldstein-Kac model \eqref{eqn:intro1} mentioned in the introduction and no waves can be formed ($\Lambda(\rho)=\rho/\lambda(\rho)$ is injective). If $\lambda(\rho)=a+b\rho^2$ is quadratic, non-trivial pairs fulfilling condition {\bf A} in Rem.~\ref{rmk:ageless} can be found as $u\,v=a/b$, however the larger of the two values will always violate condition {\bf B}.
\end{example}


\begin{example}{\em Sigmoid turning function.}
\label{exp:sigmoid_ageless} We demonstrate the construction for a sigmoid turning function for which $\Lambda(\rho)$ changes sign twice. We use $\lambda(\rho)=\lu+(\lo-\lu)\alpha \rho^2/(1+\alpha\rho^2)$ with $\lu=0.5$, $\lo=10$, $\alpha=0.125$. Note that this choice of turning function does not fulfill \eqref{eqn:symmetric}.  Fig.~\ref{fig:ageless}a shows how stable counter-propagating traveling waves (in the sense of Rem.~\ref{rmk:ageless}) are constructed: One seeks the intersections of the implicitly defined curves
$$
\Lambda(u)=\Lambda(v),\quad \Omega(u)=\Omega(v)
$$
in $(u,v)$-space (show in blue and red respectively) for which $\Lambda'>0$ (blue shading in Fig.~\ref{fig:ageless}a). To confirm the obtained values $(w_1,w_2)$ we simulated \eqref{eqn:noAge2} for various initial conditions shown in Fig.~\ref{fig:ageless}b upper row. The final traveling profiles are shown in Fig.~\ref{fig:ageless}b lower row and confirm the analytically determined stable pair $(w_1,w_2)$. Fig.~\ref{fig:ageless}c shows that for the fast system \eqref{eqn:fast} this indeed corresponds to a heteroclinic orbit connecting $w_1$ and $w_2$.
\end{example}


\begin{example}{\em Double-sigmoid turning function.}
\label{exp:doubleSigmoid} A non-trivial example where not all constructed values are actually stable is given by a double-sigmoid function, parametrized by $\lu<\hat\lambda<\lo$, $0<\underline\rho<\overline\rho$ and $0<\delta<\min\left\{\underline\rho, \frac{\overline\rho-\underline\rho}{2}\right\}$ and shown in the inset in Fig.~\ref{fig:ageless}d.
\begin{align*}
\lambda(\rho)=\begin{cases}
\lu &\rho<\underline\rho-\delta\\
\lu+\frac{\hat\lambda-\lu}{2\delta^2}\left(\rho-\underline\rho+\delta\right)^2 & \underline\rho-\delta\leq \rho<\underline\rho\\
\hat\lambda-\frac{\hat\lambda-\lu}{2\delta^2}\left(\rho-\underline\rho-\delta\right)^2 & \underline\rho\leq \rho<\overline\rho-\delta\\
\hat\lambda & \underline\rho+\delta\leq \rho<\overline\rho-\delta\\
\hat\lambda+\frac{\lo-\hat\lambda}{2\delta^2}\left(\rho-\overline\rho+\delta\right)^2 &\overline\rho-\delta\leq \rho<\overline\rho\\
\lo-\frac{\lo-\hat\lambda}{2\delta^2}\left(\rho-\overline\rho-\delta\right)^2 &\overline\rho\leq \rho<\overline\rho+\delta\\
\lo & \rho\geq \overline\rho+\delta.
\end{cases}
\end{align*}
We use as parameters $(\lu,\hat\lambda, \lo)=(2.5, 5.25,8)$, $(\underline\rho, \overline\rho)=(0.67, 1.33)$ and $\delta=1/6$.
Inspection of the curves defined by $\Omega$ and $\Lambda$ show three possible pairs of stable values, marked by 1,2,3 in Fig.~\ref{fig:ageless}d. However in light of the conserved average mass of one, pairs 2 and 3 would lead to a too low and too high overall mass respectively, hence cannot be candidates for long-term behavior. For the remaining pair 1, the orbit in phase-space shown in Fig.~\ref{fig:ageless}e shows that it does not constitute the sought after heteroclinic orbit connecting the two states. Consequently we expect to find no stable piecewise constant traveling waves for this choice of turning function, this was confirmed numerically.
\end{example}


\subsection{Traveling Waves in the Full System}
\label{ssec:waves_full}

In this section we construct traveling waves for the full system \eqref{eqn:main_scaled} for arbitrary aging and turning functions. We start by noting that Def. \ref{def:travelingWaves} does not refer to $u_1$ and $v_1$. The main result of this section is to postulate specific dependencies of $u_1$ and $v_1$ on the traveling wave frames $x\pm ct$ and be thus able to derive closed, decoupled equations for the wave profiles.


\paragraph{Construction.}
We seek counter-propagating traveling wave solutions moving with speed $\pm 1$. We expect that the composition of a wave in terms of reversible and non-reversible individuals is modified by the opposing wave. This motivates the introduction of the functions $A(x,t)=u_1(x,t)/u(x,t)$ and $B(x,t)=v_1(x,t)/v(x,t)$, i.e. the ratios of reversible individuals within a right- and left moving wave respectively. We now postulate that $A(x,t)$ depends only on $x+t$ and $B(x,t)$ only on $x-t$, i.e. the traveling wave frame of the \textit{opposing} group. The ansatz can be summarized as
\begin{align}
\label{eqn:ansatz}
&u(x,t)=P(x-t),\quad u_1(x,t)=P(x-t)A(x+t)\\
&v(x,t)=M(x+t),\quad v_1(x,t)=M(x+t)B(x-t).\nonumber
\end{align}
Substituting this into \eqref{eqn:main_scaled} and assuming that $A$ and $B$ are bounded away from zero, we find that $P(\xi)$, $M(\xi)$, $A(\eta)$ and $B(\eta)$ have to fulfill
\begin{align}
\label{eqn:totals}
\frac{P(\xi)}{\lambda(P(\xi))B(\xi)}=\frac{M(\eta)}{\lambda(M(\eta))A(\eta)}=r,
\end{align}
for some $r>0$, and
\begin{align}
\label{eqn:fractions}
 2\,A'(\eta)&=\gamma(M(\eta))\left[1-A(\eta)\right]-\lambda(M(\eta))A(\eta)\\
 -2\,B'(\xi)&=\gamma(P(\xi))\left[1-B(\xi)\right]-\lambda(P(\xi))B(\xi).\nonumber
\end{align}

\begin{figure}[t]
\centering
\includegraphics[width=\textwidth]{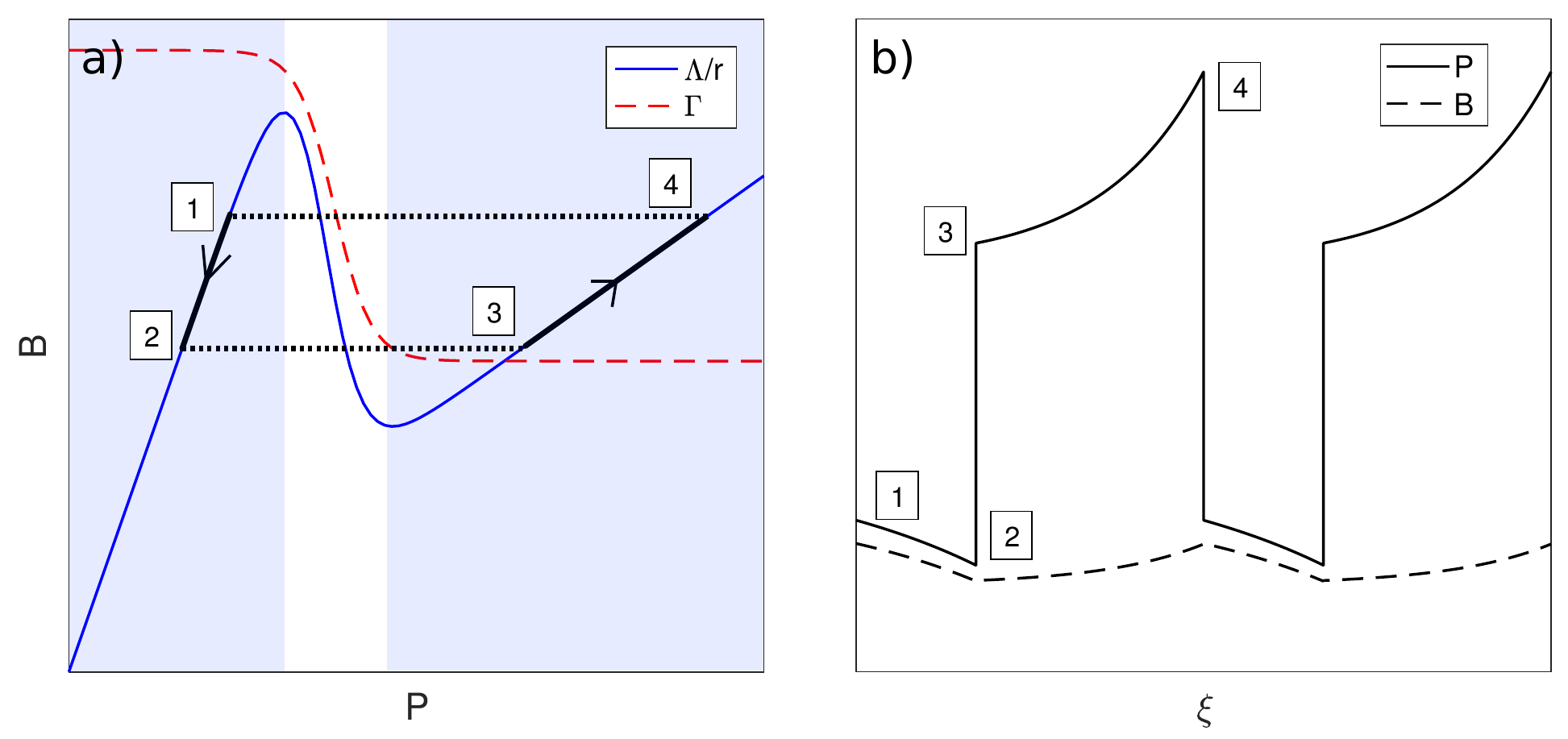}
\caption{\textit{Construction of Traveling Waves.} a) Construction of admissible solutions in the sense of Def.~\ref{def:admissible} in $(P,B)=(u,v_1/v)$ space (left) and the corresponding solutions $P(\xi), B(\xi)$ (right). Left: Solutions have to lie on $\Lambda_c^r$ (solid blue) for some fixed $r>0$. The curve $\Gamma$ (dashed red) divides the space according to the sign of $B'$. The white region marks non-admissible areas. $1\rightarrow 2$: Since $B<\Gamma$, $B'<0$ and we move downwards along $\Lambda_c^r$. $2\rightarrow 3$: We jump from one admissible branch to another one (dotted), $B$ does not change its value, but $P$ does. $3\rightarrow 4$: Now $B>\Gamma$ and thereby $B'>0$ and we move along $\Lambda_c^r$. $4\rightarrow 1$: We again jump horizontally in $(P,B)$ space back to the original branch (dotted). Right: The corresponding solutions $P(\xi)$ and $B(\xi)$, numbers refer to the corresponding numbers in $(P,B)$ space in the left figure.}
\label{fig:construction}
\end{figure}

\begin{remark}{\em Decoupling of the Wave Frames.}
Classically, when searching for traveling waves, the use of a traveling wave ansatz allows the reduction to one or several ordinary differential equations. However, for families of {\em counter-propagating} traveling waves, this procedure doesn't work, since the two wave frames would typically interact with each other. The key consequence of the above ansatz is the decoupling for the frames $\eta=x+t$ and $\xi=x-t$, thereby yielding a system of ODEs and algebraic equations. Note that this requires very particular structural properties. If, e.g. $\lambda$ or $\gamma$ were to be functions also of the other density (as suggested e.g. in \cite{Lutscher2002}), the procedure would break down.
\end{remark}
\medskip
\noindent
While the mathematical usefulness of this ansatz is obvious, its relevance will be demonstrated numerically in Sec.~\ref{sec:numerics}. In the following we will only refer to $P(\xi)$ and $B(\xi)$, since $M(\eta)$ and $A(\eta)$ can be treated analogously. If we analyze \eqref{eqn:totals} and \eqref{eqn:fractions} in the $(P,B)$ plane, we see that the dynamics can be described with the aid of two functions and their associated curves
\begin{align}
&\Lambda(\rho)=\frac{\rho}{\lambda(\rho)}, \quad \Gamma(\rho)=\frac{\gamma(\rho)}{\gamma(\rho)+\lambda(\rho)}\label{eqn:LambdaGamma}\\
&\rho\mapsto\Lambda^r_c(\rho):=\left(\rho,\Lambda(\rho)/r \right),\quad \rho\mapsto\Gamma_c(\rho):=\left(\rho,\Gamma(\rho)\right)
\end{align}
Condition \eqref{eqn:totals} means that $(P(\xi),B(\xi))$ has to lie on $\Lambda^r_c$ for all $\xi$. The dynamic on this curve is described by \eqref{eqn:fractions}, where the direction of the dynamic, i.e. the sign of $B'$ is determined by whether $(P(\xi),B(\xi))$ is above $\Gamma_c$, i.e. $B'>0$, or below, then $B'<0$. Fig \ref{fig:construction} summarizes the construction, details can be found in the proof of Prop.~\ref{prop:waveExist}.

\begin{remark}
\label{rem:ODEforP}
Wherever $P\in\mathcal{C}^2$, equations \eqref{eqn:totals} and \eqref{eqn:fractions} can be combined to one ordinary differential equation for $P$,
\begin{align}
\label{eqn:ODEforP}
P'=\frac{\lambda(P)\left[P(\lambda(P)+\gamma(P))-\lambda(P)\gamma(P)r\right]}{2\left(\lambda(P)-P\lambda'(P))\right)}.
\end{align}
However, as noted in Prop.~\ref{prop:waveExist}, $P$ will in general not be continuous everywhere and the construction using $B$ as described in the caption of Fig.~\ref{fig:construction} allows for more insights.
\end{remark}

We will see in Sec.~\ref{sec:numerics} that the constructed waves can be observed as long-term limits of the simulated densities, indicating that at least some are stable. In this work we do not deal with the stability of the traveling waves, however motivated by the stability results in Sec.\ref{sec:steady} and the observed numerical behavior, we only focus on one class of waves, which we will call \textit{admissible} (see Def.~\ref{def:admissible}). We expect that all stable waves will be admissible, however a proof (or dis-proof) will be the subject of future work. Before stating the existence result, we therefore introduce the following useful definitions.


\begin{definition}\label{def:admissible} Admissibility and Reachability. Let $\rho>0$, $I\subset [0,\infty)$ and $\Lambda$ given by \eqref{eqn:LambdaGamma}.
\begin{enumerate}
\item We call $\rho$ \textbf{admissible}, if $\Lambda'(\rho)>0$ ($\iff \lambda'(\rho)\rho<\lambda(\rho)$).
\item We define the \textbf{reachability set} of $\rho$, $\mathcal{R}_\rho$ and that of $I$, $\mathcal{R}_I$ by
\begin{align*}
\mathcal{R}_\rho:=\{\rho^*>0\, |\, \rho^*\neq \rho,\, \Lambda(\rho)=\Lambda(\rho^*), \, \rho^*\, \text{admissible}\},\quad  \mathcal{R}_I:=\bigcup_{\rho\in I}\mathcal{R}_\rho.
\end{align*}
If $\mathcal{R}_\rho\neq \emptyset$ we call $\rho$ \textbf{reachable}.
\item We call a pair of functions $0< P(\xi)$ and $0< B(\xi)\leq 1$ \textbf{admissible solution} if there exists an $r>0$ such that they solve \eqref{eqn:totals}, \eqref{eqn:fractions} and there exists $\bar\xi>0$ such that $B(0)=B(\bar\xi)$, $P(0)=P(\bar\xi)$ and $P(\xi)$ is admissible for all $\xi\in[0,\bar\xi]$.
\end{enumerate}
\end{definition}

The following proposition shows the existence of admissible, counter-propagating waves and details several properties. The most notable is that $P$ is necessarily non-continuous.

\begin{proposition}
\label{prop:waveExist} Let $\gamma, \lambda \in \mathcal{C}^k([0,\infty))$ $k\geq 2$ and assume that there exist constants $\lu, \underline\gamma>0$, such that 
$$\underline \gamma\leq \gamma(\rho),\quad \lu \leq \lambda(\rho),\quad \forall \rho\geq 0.$$
Further assume that there exist values $0<\rho_1<\rho_2<\rho_3$ such that
\begin{align*}
\begin{cases}
&\Lambda'(\rho)>0\quad  \rho\in(0,\rho_1)=:I_1 \\\
&\Lambda'(\rho) <0\quad  \rho\in(\rho_1,\rho_2)\\
&\Lambda'(\rho) >0\quad  \rho\in(\rho_2,\rho_3)=:I_2.
\end{cases}
\end{align*}
There there exists an admissible solution in the sense of Def.~\ref{def:admissible}, $B\in\mathcal{C}([0,\infty))$ and piecewise in $\mathcal{C}^k([0,\infty))$ and $P$ piecewise in $\mathcal{C}^k([0,\infty))$, but $P\notin \mathcal{C}([0,\infty))$. For all values $\hat\xi\in[0,\infty)$ where $B$ and $P$ are not in $\mathcal{C}^k$, $P$ fulfills the jump condition
\begin{align*}
\lim_{\xi\rightarrow\hat\xi^+}\frac{P(\xi)}{\lambda(P(\xi))}=\lim_{\xi\rightarrow\hat\xi^-}\frac{P(\xi)}{\lambda(P(\xi))}.
\end{align*}
\end{proposition}


\begin{figure}[t]
\centering
\includegraphics[width=0.7\textwidth]{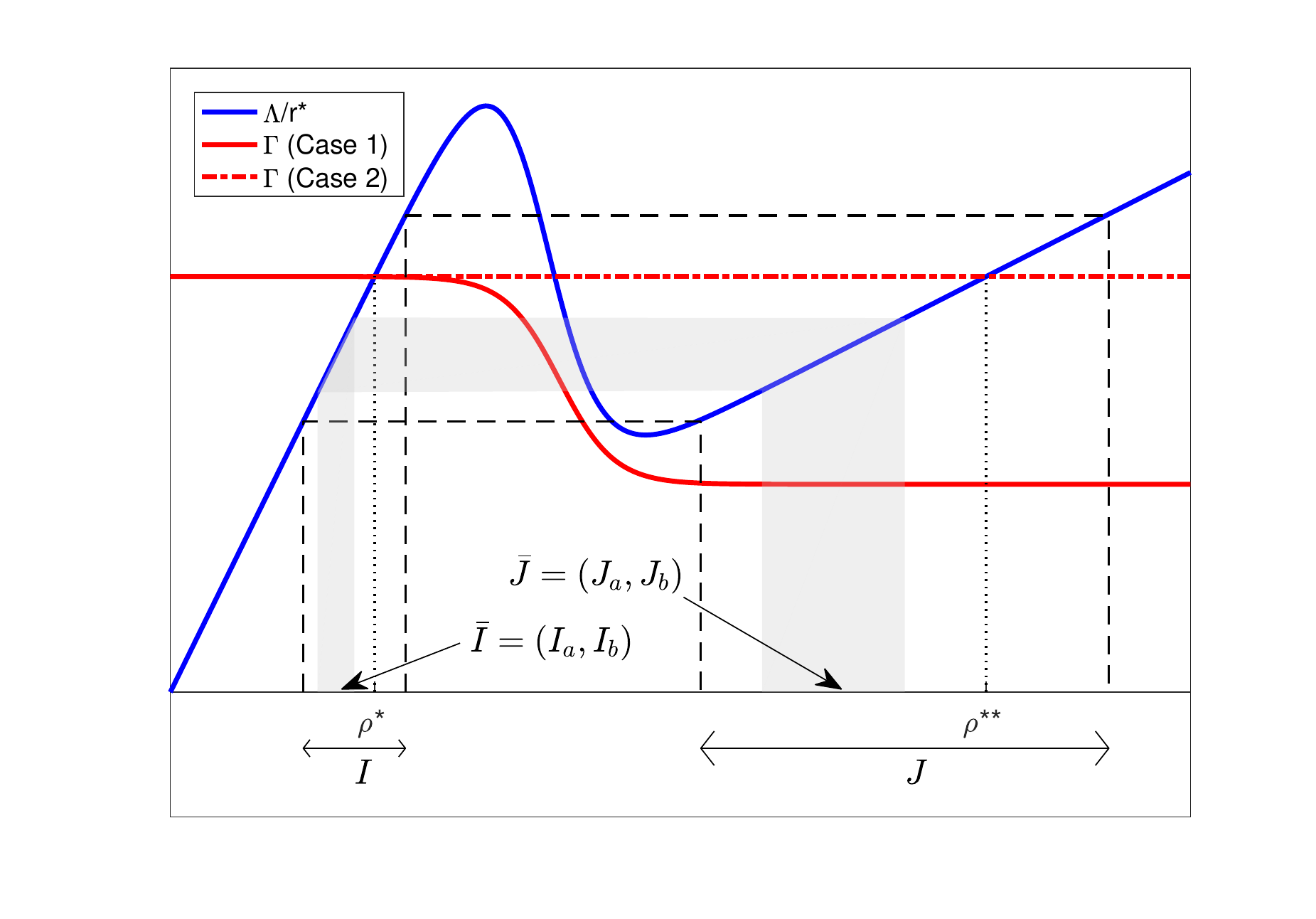}
\caption{\textit{Proof of Existence.} }
\label{fig:proofExist}
\end{figure}

\medskip 
\noindent 
{\bf Proof.} Please refer to Fig.~\ref{fig:proofExist} for an illustration of the the proof. Since $\Lambda(0)=0$ and $\Lambda$ makes at least two changes in direction, we can find open intervals $I=(a,b)\subset I_1$ and $J\subset I_2$, such that
$$
\Lambda(I)=\Lambda(J)
$$
We choose $I$ and $J$ such that their interval ends do not include local extrema of $\Lambda$. We now search for admissible solutions, with $P(\xi)\in I\cup J$ for all $\xi>0$. We claim that it is always possible to find $r>0$, such that
\begin{equation}
\label{eqn:intersect}
\Gamma(\rho)=\Lambda(\rho)/r
\end{equation}
has a solution $\rho\in I$. This can be seen easily by noting that if we define 
$$r_1=\frac{\Lambda(a)}{\Gamma(a)},\quad r_2=\frac{\Lambda(b)}{\Gamma(b)}$$
the intermediate value theorem applied to $\rho\mapsto \frac{\Gamma(\rho)}{\Lambda(\rho)}$ together with the continuity of $\gamma$ and $\lambda$ ensures that \eqref{eqn:intersect} has a solution $\rho\in (a,b)$ for each
\begin{align*}
r\in(\min\{r_1,r_2\}, \max\{r_1,r_2\})=:R_1.
\end{align*}
We require that for admissible solutions $B\leq 1$. In light of \eqref{eqn:totals}, we therefore require $r\in(\Lambda(b), \infty)=:R_2$. Since $\Lambda(b)<r_2$, 
$$R_1\cap R_2=:R \neq \emptyset$$
and we denote by $r^*$ an element in $R$ and define $\rho^*$ by $\Gamma(\rho^*)=\Lambda(\rho^*)/r^*$. W.l.o.g. we assume that $\Gamma$ crosses $\Lambda/r^*$ at $\rho^*$ (as opposed to just touching it), hence the sign of $B'$ in \eqref{eqn:fractions} changes at $P=\rho^*$. Due to our construction of $I$ and $J$ there exists a $\rho^{**}\in J$ with $\Lambda(\rho^*)=\Lambda(\rho^{**})$.
\vspace{0.2cm}
\newline
\textit{Case 1:} If $\Gamma(\rho^{**})\neq \Lambda(\rho^{**})/r^*$, we can find subintervals $\bar I =(I_a, I_b) \subset I$ and $\bar J= (J_a, J_b)\subset J$ such that
$$
\Lambda(\bar I)/r^*=\Lambda(\bar J)/r^*
$$
and $B'$ in \eqref{eqn:fractions} has opposite signs for $P\in \bar I$ and $P\in \bar J$. We assume $B'<0$  in $\bar I$ and $B'>0$ in $\bar J$. In light of Rmk.~\ref{rem:ODEforP} we can then define an admissible solution by
\begin{enumerate}
\item Set $P(0)=I_b$, solve \eqref{eqn:ODEforP} with $r=r^*$ until $P(\xi_1)=I_a$ for some $\xi_1>0$.
\item Set $P(\xi_1)=J_a$ and solve \eqref{eqn:ODEforP} with $r=r^*$ until $P(\bar\xi)=J_b$  for some $\bar \xi>\xi_1>0$.
\item Define $B(\xi)=\Lambda(P(\xi))/r^*$.
\end{enumerate}
This procedure works since on $\bar I$ and $\bar J$ we have $\lambda(\rho)>\rho\lambda'(\rho)$ and $\sign P'=\sign B'$. Further we can bound
\begin{align*}
&P'\leq -\hat P_{\bar I}<0\,\,\text{on $\bar I$},\quad P'>\hat P_{\bar J}>0\,\,\text{on $\bar J$},\quad \text{where}\\
&\hat P_W=\frac{\lu^2(\lu+\underline\gamma)r^*}{2\sup_{Q\in W}(\lambda(Q)-Q\lambda'(Q))}\inf_{W}|\Lambda/r^*-\Gamma|>0,
\end{align*}
hence we will be able to reach $I_a$ and $J_b$ for finite $\xi_1$ and $\bar\xi$. Whenever $P$ doesn't jump it is the classical solution of an ODE with the right hand side in $\mathcal{C}^{k-1}$, hence away from discontinuities of P, we have $P,B\in \mathcal{C}^k$. Further since $\Lambda(I_a)=\Lambda(J_a)$ and $\Lambda(I_b)=\Lambda(J_b)$, the resulting $B$ will be continuous.
\vspace{0.2cm}
\newline
\textit{Case 2:} If $\Gamma(\rho^{**})= \Lambda(\rho^{**})/r^*$ we define
$$
B(\xi)\equiv \Lambda(\rho^*)/r^*, \quad P\,\, \text{piecewise constant with}\,\, P(\xi)\in\{\rho^*, \rho^{**}\}.
$$
Since $\Lambda(\rho^*)/r^*$ is a fixpoint of \eqref{eqn:fractions}, $B$ fulfills \eqref{eqn:fractions} trivially and we have defined an admissible solution.
\qed



\begin{lemma}
\label{lem:bounds}
Let the assumptions of Prop.~\ref{prop:waveExist} hold and let $P,B$ be admissible solutions in the sense of Def.~\ref{def:admissible}. Further let
\begin{align}
&\overline P:=\max \left\{ \{\rho_M\}\cup \mathcal{R}_{\rho_M}|\, \Lambda'(\rho_M)=0,\,  \Lambda''(\rho_M)<0,\, \rho_M\,\, \text{is reachable}\right\},\\
&\underline P:=\min \left\{ \{\rho_M\}\cup \mathcal{R}_{\rho_M}|\, \Lambda'(\rho_M)=0,\,  \Lambda''(\rho_M)>0,\, \rho_M\,\, \text{is reachable}\right\}.
\end{align}
Then we have 
$$
\underline P\leq ||P||_\infty\leq \overline P.
$$
For the fraction $B$ we can conclude that
$$
\inf_{\rho\in[0,\infty)}\Gamma(\rho)\leq||B||_\infty\leq \sup_{\rho\in[0,\infty)}\Gamma(\rho)<1.
$$
\end{lemma}
\medskip 
\noindent 
{\bf Proof.} The proof consists of simple geometrical arguments using the construction depicted in Fig.~\ref{fig:construction}. The bounds of $P$ result from examining $\Lambda(\rho)$ in terms of reachability, disregarding $\Gamma(\rho)$. The borders of reachable intervals are determined by the minima and maxima of $\Lambda(\rho)$. The bounds of $B$ are a result of the observation that the solution curve in $(P,B)$ has to be able to jump to a branch of $\Lambda(\rho)/r$ with opposite sign of $B'$. \qed


\begin{example}{\em A simple sigmoidal $\lambda$}. To apply the construction to a simple case we set $\gamma(\rho)\equiv\gamma$ and use the simplest continuous approximation of a sigmoid $\lambda(\rho)$ given by
\begin{equation}
\label{eqn:lambda_piecewise}
\lambda(\rho)=\begin{cases}
\underline\lambda\quad & \rho<1-\e \\
\underline\lambda+\frac{\rho-1+\e}{2\e}(\overline\lambda-\underline\lambda) \quad & 1-\e\leq \rho < 1+\e\\
\overline\lambda\quad & \text{else}.
\end{cases}
\end{equation}
Here $\e>0$ is a small parameter, for $\e\rightarrow 0$ $\lambda(\rho)$ would converge to a piecewise constant function. Strictly speaking $\lambda$ doesn't fulfill the smoothness requirements of Prop.~\ref{prop:waveExist}, however the result could easily be generalized to continuous functions, so we won't let that stop us. The instability condition \eqref{eqn:stability} requires
$$
\e<\frac{\overline\lambda-\underline\lambda}{\overline\lambda+\underline\lambda}
$$
to allow for admissible counter-propagating traveling waves. Application of Lem.~\ref{lem:bounds} gives
$$
\frac{\underline\lambda}{\overline\lambda}(1+\e)<P<\frac{\overline\lambda}{\underline\lambda}(1-\e),\quad \frac{\gamma}{\gamma+
\overline\lambda}<B<\frac{\gamma}{\gamma+\underline\lambda},
$$
Further, referring to the constant $r$ in \eqref{eqn:totals} and Prop. \ref{prop:waveExist}, we can conclude that
$$
\quad \frac{\gamma+\underline\lambda}{\gamma\overline\lambda}(1-\e)<r< \frac{\gamma+\overline\lambda}{\gamma\underline\lambda}(1+\e).
$$

\paragraph{Shapes.} We apply the construction in Prop.~\ref{prop:waveExist} using \eqref{eqn:ODEforP} to obtain the following wave crest and trough  shape: Let $P_{c,0}>r\frac{\gamma\overline\lambda}{\gamma+\overline\lambda}$ and $P_{t,0}<r\frac{\gamma\underline\lambda}{\gamma+\underline\lambda}$ be the initial densities of the crest and trough respectively, then
\begin{align*}
&P_\text{crest}(\xi; P_{c,0})=r\frac{\gamma\overline\lambda}{\gamma+\overline\lambda}+e^{\frac{\gamma+\overline\lambda}{2}\xi}\left(P_{c,0}-r\frac{\gamma\overline\lambda}{\gamma+\overline\lambda}\right)\\
&P_\text{trough}(\xi; P_{t,0})=r\frac{\gamma\underline\lambda}{\gamma+\underline\lambda}+e^{\frac{\gamma+\underline\lambda}{2}\xi}\left(P_{t,0}-r\frac{\gamma\underline\lambda}{\gamma+\underline\lambda}\right).
\end{align*}
In particular wave crest are convex and wave troughs are concave. When connecting wave crests with troughs, the continuity of $B$ leads to a jump conditions for $P$ at the discontinuity $\xi_c$
\begin{align*}
& \lim_{\xi\rightarrow \xi_c-}P(\xi)=\frac{\overline\lambda}{\underline\lambda}\lim_{\xi\rightarrow \xi_c+}P(\xi)\quad \text{from crest to trough}\\
& \lim_{\xi\rightarrow \xi_c-}P(\xi)=\frac{\underline\lambda}{\overline\lambda}\lim_{\xi\rightarrow \xi_c+}P(\xi)\quad \text{from trough to crest}
\end{align*}
In particular we get for a wave with a crest for $\xi\in[0,\xi_1]$, then followed by a trough for $\xi\in[\xi_1, \xi_2]$
\begin{align*}
&P(\xi)=\begin{cases}
P_\text{crest}(\xi; P_{c,0}), \quad 0<\xi<\xi_1\\
P_\text{trough}(\xi-\xi_1; P_{t,0}), \quad \xi_1<\xi<\xi_2
\end{cases},\quad \text{where}\quad P_{t,0}=P_\text{crest}(\xi_1,P_{c,0})\frac{\underline\lambda}{\overline\lambda}.
\end{align*}
Requiring the wave to be periodic, i.e. $B(\xi_2)=B(0)$ or equivalently $P(\xi_2)=P_{c,0}\frac{\underline\lambda}{\overline\lambda}$ yields an expression for $P_{c,0}$ in terms of $\xi_1$ and $\xi_2$. If one additionally fixes the mass of the wave, also $r$ is determined.
\end{example}


\begin{figure}[t]
\centering
\includegraphics[width=\textwidth]{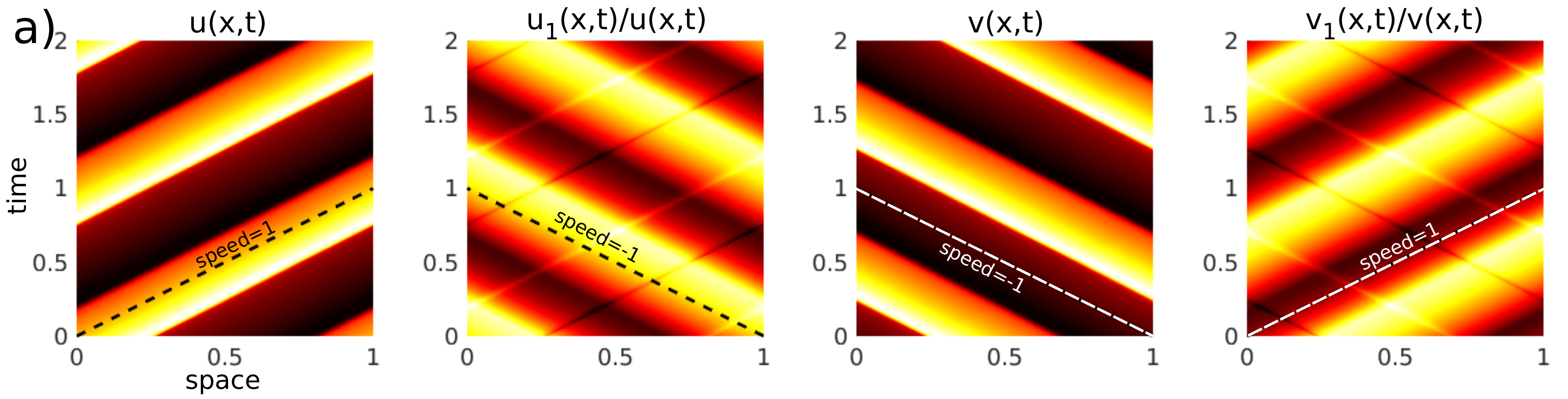}
\caption{\textit{Travling Wave Frames.} a) Space-time plots for the densities $(u,u_1/u, v, v_1/v)=(P,A,M,B)$ for a simulation using the sigmoid turning from Ex.~\ref{exp:Hopf}. As initial conditions the isotropic steady state was perturbed with a sine-wave. Shown is the result for the time $t\in[48,50]$. Colors represent high (white) and low (black) densities. Dashed lines make speed$=\pm 1$.}
\label{fig:spaceTime}
\end{figure}

\begin{figure}[p]
\centering
\includegraphics[width=\textwidth]{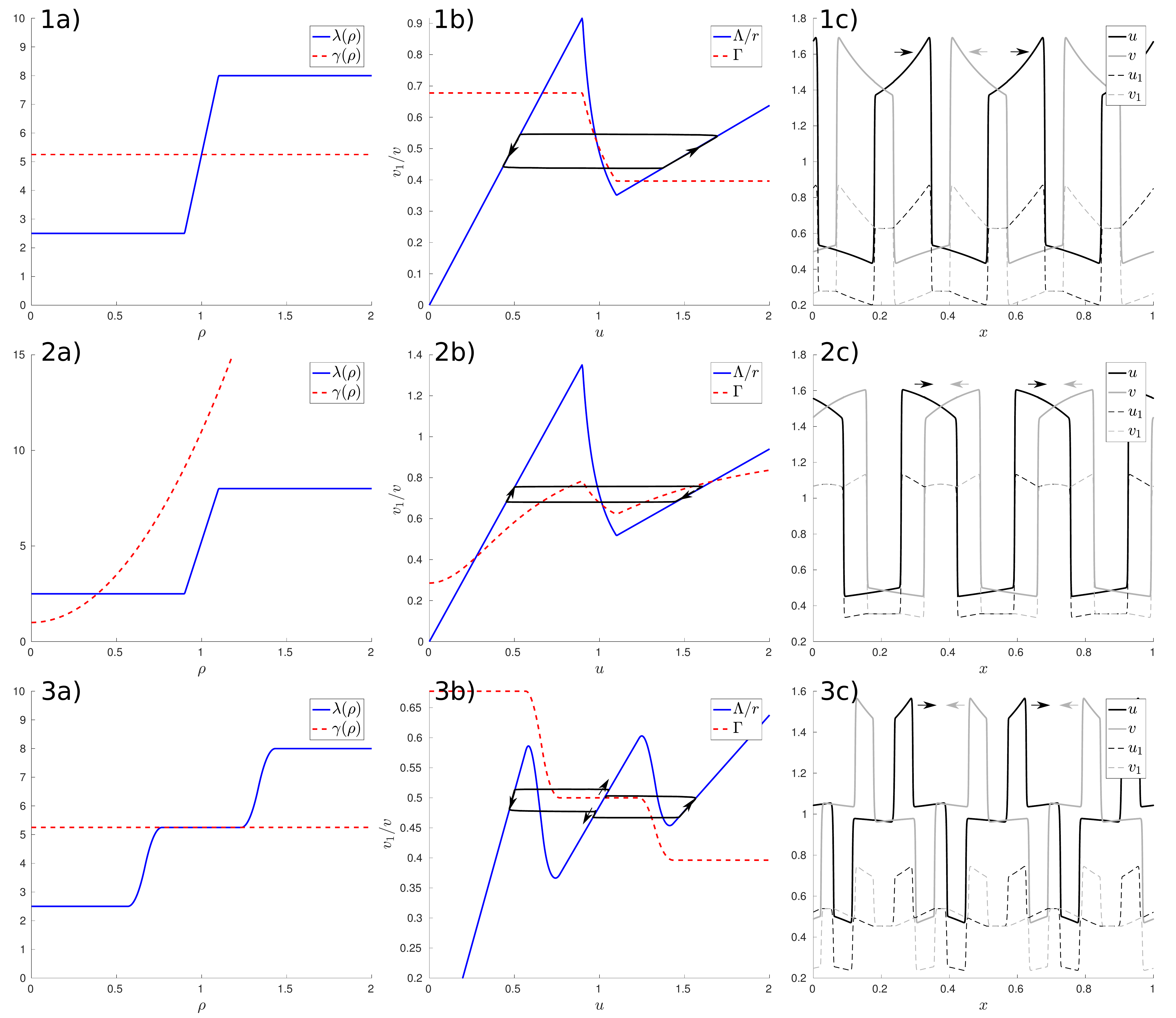}
\caption{\textit{Wave Gallery.} Depicted are three example choices of $\lambda(\rho)$ and $\gamma(\rho)$ (first column, letter a), the corresponding wave construction in $(P,B)=(u, v_1/v)$-space (middle column, letter b, arrows mark directions) and the resulting simulated waves (last column, letter c). Example 1: Constant $\gamma(\rho)$ and piecewise linear $\lambda(\rho)$. Wave crests increase in travel direction, wave troughs decrease. Example 2: Quadratic $\gamma(\rho)$ and piecewise linear $\lambda(\rho)$. Now wave crests decrease in travel direction, wave troughs increase. Example 3: Constant $\gamma(\rho)$ and $\lambda(\rho)$ having a triple-step shape. Wave shapes are more complicated with four piecewise smooth segments. As initial conditions in all cases I used the isotropic steady state plus a sine/cosine wave.}
\label{fig:gallery}
\end{figure}

\section{Simulations}
\label{sec:numerics}

We simulated model \eqref{eqn:main_scaled} in one space dimension with periodic boundary conditions. The transport and reaction terms were implemented using operator splitting with explicit upwind or downwind (for the right- and left moving family respectively) finite differences for the transport term and an explicit treatment of the reaction term. We also tested implementing the reaction term with an explicit Runge-Kutta (4,5) formula using the ode45 solver of Matlab and a Lax-Friedrichs Scheme for the transport term, both with no significant gains in accuracy. We used a spatial step of $\Delta x=6.25\times 10^{-4}$ and a time step of $\Delta t=6.1875\times 10^{-4}$ which leads to only a very small amount of numerical diffusion introduced through the discretization of the transport term.

\paragraph{Confirmation of Traveling Wave Frames.} The crucial assumption that allowed for the construction of the traveling waves for system \eqref{eqn:main_scaled} was that both $u$ and $v_1/v$ form traveling waves with speed 1 to the right (and the left for $v$ and $u_1/u$). We demonstrate the validity of this ansatz by examining the space-time evolution for the sigmoid turning function from Ex.~\ref{exp:Hopf}. As initial conditions we use the isotropic steady state solutions \eqref{eqn:SS1} perturbed by a sine function with period 1. Fig.~\ref{fig:spaceTime} confirms several assumptions and results from Sec.~\ref{ssec:waves_full}: The fractions $u_1/u$  and $v_1/v$ form traveling waves moving in the direction {\em opposite} of their transport direction and $u$ and $v$ traveling waves moving in their transport direction. All observed speeds are $\pm 1$, as postulated. We further see that $u$ and $v$ appear to be piecewise continuous, whereas $u_1/u$ and $v_1/v$ appear continuous.

\paragraph{Wave Gallery.} Fig.~\ref{fig:gallery} shows a gallery of different wave shapes produced by the dynamics for three choices of aging and turning functions. For each case the correct value for $r$ in \eqref{eqn:totals} was determined numerically. We see that the construction gives correct wave shapes even for complicated choices for $\gamma(\rho)$ and $\lambda(\rho)$ and that the wave profiles include growing crests (upper row), falling crests (middle row) and combinations of both (lower row).

\section{Conclusion}

In this work we have discussed a simple system of four transport-reaction equations, coupled by their nonlinear reaction terms. The system was derived in \cite{DegMan2017} as an application to the rippling behavior of myxobacteria, but should be interpreted in the larger context of pattern formation in hyperbolic systems (\cite{primi2009,hadeler1999}). Despite its simplicity the system's behavior is very rich and we showed that it includes both pulsating-in-time, constant-in-space solutions as well as highly non-trivial counter-propagating traveling waves, whose explicit construction was detailed. The key to uncoupling the two wave frames is the intriguing feature, that the fractions of non-reversing/reversing right-moving densities form traveling waves moving against their transport direction.
\newline\par
Several mathematically and biologically interesting questions arise from this study. The first one concerns the stability of the counter-propagating traveling waves: In the limit of fast-aging two conditions determine the selected wave heights. 1) Linear stability of the individual heights led to an inequality and 2) regularization with a diffusion term yielded an integral condition, which was necessary to fully characterize the wave heights selected in simulations. For the full system \eqref{eqn:main_scaled} both the inequality and the integral condition remain to be found. One expects the inequality condition to be similar or equal to the admissibility condition in Def.~\ref{def:admissible}. The numerical results suggest that also an analogue of the the integral condition exists for the full system. Secondly, at this point results about wavelength selection are missing. In \cite{DegMan2017} it was noted that a diffusion term might be necessary to avoid very short-length waves. As pointed out in \cite{igoshin2004,primi2009} linear stability analysis will not be the correct tool to predict wavelengths for such systems, alternative mechanisms for wavelength selection in the memory-free model were discussed in \cite{scheel2017}, but it remains to be seen if they can be generalized to the 2-age system.  Finally it will be interesting to explore if a similar ansatz to decouple waveframes can be used for the continuous-age model \eqref{eqn:intro3}.

\section{Acknowledgments.}

The author wants to express gratitude to P. Degond and C. Schmeiser for the insightful discussions and, A. Mogilner for his support.

\bibliographystyle{abbrv}

\end{document}